%




\magnification=\magstep1
\documentstyle{amsppt}
\def\CC{{\Bbb C}}

\font\sc=cmcsc10

\def\ds{\displaystyle}

\def\delbar{{\bold{\bar\partial}}}
\def\Dbar{{\,\bbar{\italic{\!D}}}}
\def\Vbar{{\,\bbar{\italic{\!V}}}}
\def\cinfty {$C^\infty$}
\def\zoform {$(0, 1)$-form}
\def\zoforms {$(0, 1)$-forms}

\def\lpee { \l_p }
\def\RR {{\Bbb R}}
\def\ZZ {{\Bbb Z}}
\def\KK {{\Cal K}}
\def\OO {{\Cal O}}
\def\LL {{\Cal L}}
\def\gg {{\frak g}}
\def\delbarclosed{$\delbar$-closed}
\def\delbarexact{$\delbar$-exact}
\def\dbe{$(1.1)$\ }
\def\ro{\varrho}

\def\epsz{\varepsilon}
\def\fii{\varphi}
\def\fv{func\-tion}
\def\fvs{func\-tions}
\def\holo{hol\-o\-mor\-phic}
\def\homog{ho\-mo\-ge\-ne\-ous}
\def\mhomog{multi\-ho\-mo\-ge\-ne\-ous}
\def\mh{\homog}
\def\mind{multi\-index}
\def\kmh{$k$-\mh}

\def\cohom{co\-hom\-o\-logy}
\def\gp{group}
\def\isom{isomor\-phism}

\def\thm{theorem}
\def\acc{almost complex}
\def\Acc{Almost complex}
\def\mfold{manifold}
\def\cpx{complex}
\def\struc{struc\-ture}
\def\nbd{neighbor\-hood}
\def\halmaz#1 {{\left\{#1\right\}}}
\def\nm#1 {{\left\Vert #1\right\Vert}}
\def\bspc{Banach space}
\def\seq{sequence}
\def\dim{dimen\-sion}
\def\cvne{C^\infty_{0,1}}

\def\la{\lambda}
\def\sa{\sigma}
\def\WLOG{Without loss of generality}

\def\bbar#1{{\bold{\bar{#1}}}}
\def\za{\zeta}
\def\en{{1,0}}
\def\ne{{0,1}}
\def\Tne{{T^\ne}}
\def\Vne{{V^\ne}}
\def\Hne{{H^\ne}}
\def\VWne{{V^\ne_W}}
\def\HWne{{H^\ne_W}}

\def\ctm{{\CC\otimes TM}}
\def\ek{{\textstyle{1\over 2}}}

\def\acm{{3.1}}
\def\Prop{Proposition}
\def\hivI{Theorem~1.1}
\def\hivII{Proposition~1.2}
\def\hivIII{Theorem~2.1}
\def\hivIV{Proposition~2.2}
\def\hivV{Theorem~2.3}
\def\hivVI{Proposition~2.4}
\def\hivVII{\Prop~2.5}
\def\hivVIII{\Prop~2.6}
\def\hivIX{\Prop~2.7}
\def\hivX{Theorem~3.1}
\def\hivXI{\Prop~3.2}

\pageheight{9truein}
\topmatter
\title
	On the $\delbar$-equation in a Banach space 
\endtitle
\comment
	Sur l'\'equation $\delbar$ dans un espace de Banach
\endcomment

\author
	Imre Patyi\footnote""{This research was partially supported
	by an NSF grant.\hfill\hfill}
\endauthor
\address{Department of Mathematics, Purdue University,
	 West Lafayette, IN 47907--1395, USA}
\endaddress
\abstract
{
	We define a separable Banach space  $X$  and prove the existence
	of a $\delbar$-closed \cinfty-smooth \zoform\  $f$  on the
	unit ball  $B$  of  $X$,  which is not  $\delbar$-exact on 
	any open subset.
	Further, we show that the sheaf \cohom\ \gp{}s
	$H^q(\Omega, \OO)=0$,  $q\ge1$, where  $\OO$  is
	the sheaf of germs of \holo\ \fvs\ on  $X$,  and
	$\Omega$ is any pseudoconvex domain in  $X$, e.g.,
	$\Omega=B$.
	As the Dolbeault \gp\ $H^{0,1}_{\;\delbar}(B)\not=0$,
	the Dolbeault \isom\ \thm\ does not generalize to
	arbitrary \bspc{}s.	
	Lastly, we construct a \cinfty-smooth integrable
	\acc\ \struc\ on $M=B\times\CC$
	such that no open subset of  $M$  is
	bi\holo\ to an open subset of a \bspc.
	Hence the Newlander--Nirenberg theorem
	does not generalize to arbitrary Banach \mfold{}s.

\comment
	On d\'efinit un espace de Banach s\'eparable $X$
	et on montre l'existence d'une forme
	$\delbar$ ferm\'ee du type $(0,1)$ de classe \cinfty\
	sur la boule--unit\'e $B$ de $X$, qui n'est pas $\delbar$ exacte
	dans aucun ouvert.
	  Au m\^eme temps, $H^q(\Omega, \OO)=0$, si $q\ge1$ et
	$\Omega\subset X$ ouvert pseudo-convexe, p.~ex., $\Omega=B$.
	Il en suit que l'ismorphisme de Dolbeault ne se g\'en\'eralise pas
	aux espace de Banach quelconques.
	  On montre \'egalement que le th\'eor\`eme de 
	Newlander--Nirenberg ne se g\'en\'eralise pas aux vari\'et\'es de
	Banach quelconques.
\endcomment

}
\endabstract

\endtopmatter
{\hfill\it \'Edesany\'amnak, \'Edesap\'amnak.}
\bigskip

\TagsOnLeft
\baselineskip=12pt
\parskip=4pt

\subhead Introduction
\endsubhead

	This paper addresses three fundamental problems that arise
	in complex analysis on \bspc{}s and on Banach manifolds.

	The first concerns vanishing of Dolbeault \cohom\ \gp{}s.
	Presently there is one definitive result on this:
	the Dolbeault \cohom\ \gp\ $H^{0,1}_{\;\delbar}(\Omega)=0$
	for any pseudoconvex open $\Omega\subset\l_1$,
	see [L3, Corollary~0.2]. In no other infinite dimensional
	\bspc\ is a similar result available. Here, we shall show 
	that such a vanishing theorem cannot be true in complete
	generality. In Section~1 we shall define a separable	
	\bspc\ $X$ and a \cinfty-smooth \delbarclosed\ \zoform\
	$f$ on its unit ball such that on no open set $G$ is
	the equation $\delbar u=f|_G$ solvable. This implies	
	$H^{0,1}_{\;\delbar}(\Omega)\not=0$ for any bounded open
	set $\Omega\subset X$. We note that globally
	non-solvable $\delbar$-equations in Fr\'echet spaces
	were constructed earlier by Dineen~[D] and Meise--Vogt~[MV].

	The second issue to be considered is that of an 
	infinite dimensional version of the Dolbeault isomorphism 
	between the Dolbeault \cohom\ \gp{}s 
	$H^{0,q}_{\;\delbar}(\Omega)$ and the sheaf \cohom\ \gp{}s
	$H^q(\Omega, \OO)$, where $\OO$ is the sheaf of germs of \holo\
	\fvs\ on $X$.  Currently no instance of such an isomorphism
	is known when $\Omega$ is open in an infinite dimensional
	\bspc. We shall show that $H^q(\Omega, \OO)=0$, $q\ge1$,
	for all pseudoconvex open subsets $\Omega$ of the above space $X$.
	In particular, 
	$0=H^1(\Omega, \OO)\not\equiv H^{0,1}_{\;\delbar}(\Omega)\not=0$
	for any bounded pseudoconvex open set $\Omega\subset X$.
	The vanishing of sheaf \cohom\ follows from a theorem of Lempert
	[L3, Theorem~0.3] plus a Runge--type approximation theorem
	to be proved in Section~2.

	The last issue to be addressed concerns the extension of the 
	Newlander--Nirenberg theorem on integrating \acc\ \struc{}s
	to an infinite dimensional setting.
	The question is whether a (\cinfty--smooth) formally integrable
	\acc\ \mfold\ is locally bi\holo\ to a vector space.
	In finite dimensions it is true, see [NN], while in some 
	Fr\'echet manifolds it is known to be false, see [LB, L5].
	This failure in itself is perhaps not surprising, as on
	Fr\'echet manifolds even real vector fields may not be integrable.
	However, in Section~3, given any \cinfty-smooth \delbarclosed\
	but nowhere \delbarexact\ \zoform\ $f$ on the unit ball  $B$
	of a \bspc\  $X$, we explicitly construct a \cinfty-smooth
	integrable \acc\ \struc\ on  $M=B\times\CC$
	such that no open subset of $M$ is bi\holo\ to
	an open subset of a \bspc, giving a Banach \mfold\
	to which the Newlander--Nirenberg \thm\ does not
	generalize. The \mfold\ $M$ is
	a smoothly trivial principal $(\CC,+)$ bundle over $B$
	and its \acc\ \struc\ will be determined by the form $f$,
	which we use as a
	deformation tensor. 

	Below we shall use freely basic notions of infinite dimensional
	complex analysis, see [L1, L2] for the definition and 
	basic properties of the following items:
	differential calculus on infinite \dim{}al spaces; 
	smoothness classes  $C^m(\Omega)$, $C^{\,m}_{p,q}(\Omega)$
	of functions and of $(p,q)$-forms with $m=0,1,\ldots,\infty,\omega$;
	the $\delbar$-complex; complex manifolds, \acc\ \mfold{}s;
	pseudoconvexity;
	\holo\ mappings and integrability of \acc\ \struc{}s.
	
\subhead
 	Notation
\endsubhead

	Denote by  $B_X(a,r)=\halmaz{ x\in X: \nm {x-a} <r } $ the open ball
	with center  $a\in X$  of radius  $0<r\le\infty$ in a Banach space
	$(X,\nm{\cdot} )$. Put $B_X(r)=B_X(0,r)$.
	Denote by  $C^m(\Omega)$, $C^{\,m}_{0,1}(\Omega)$
	the space of complex functions and of \zoforms\ of
	smoothness class $m=0,1,\ldots,\infty,\omega$, and
	define for $u\in C^m(\Omega)$, $m<\infty$,
	the  $C^m(\Omega)$  norm by
$$
	\nm{u} _{C^m(\Omega)}=\sum_{k\le m}\sup_{x\in\Omega}\nm {u^{(k)}(x)}
	\le\infty
 $$
	where $\nm {u^{(k)}(x)} $ is the operator norm of the 
	$k$th Fr\'echet derivative
	$u^{(k)}$ of $u$.  The $C^m(\Omega)$ norm of 
	$f\in C^{\,m}_{0,1}(\Omega)$ is
	defined by
$$
	\nm{f} _{C^m(\Omega)}=\nm{u} _{C^m(\Omega\times B_X(1))}\le\infty
 $$
 	where $u(x, \xi)=f(x)\xi$ for $x\in\Omega$, $\xi\in B_X(1)$.

\subhead
	1.\ \ Non-solvability
\endsubhead

	We consider the solvability of
$$
	\delbar u = f\hbox{\ \  on\ }\Omega,\tag1.1
 $$
	where  $f\in C^{\,m}_{0,1}(\Omega)$  
	is a \delbarclosed\  \zoform\ with 
	$m=1, 2, \dots, \infty$
	on a domain  $\Omega$ in a Banach space  $X$.

	Coeur\'e in [C] (see also [M]) gave an  $f$  on $X=\Omega=l_2$
	of class  $C^1$  for which \dbe is not solvable
	on any open set.  Lempert in [L2] extended
	Coeur\'e's example and produced, with $p=2, 3, \dots$,
	a \delbarclosed\ form  $f\in C^{p-1}_{\,0,\,1}(\lpee)$ 
	for which \dbe is not solvable
	on any open set.  Based on the mere existence of
	these examples, we prove that there is a form  $f$
	of class \cinfty\ on  $\Omega=B_X(1)$
	in, say, the $l_1$-sum  $X$  of a suitable sequence of
	$\lpee(\CC^{\,n(p)})$
	spaces
	with  $p\ge 2$ integer, for which
	\dbe is not solvable on any open subset of  $\Omega$.

	Let  $Y$  be  $l_q$,  $1\le q <\infty$, or $c_0$.
	We define the  $Y$-sum  $X$  of a \seq\ of \bspc{}s
        $(X_n,\nm{\cdot} _n)_{n=1}^\infty$ by
$$
	X=\halmaz{x=(x_n):x_n\in X_n, y=|x|\in Y, \nm x =\nm y _Y} ,
 $$
	where  $|x|=(\nm {x_1} _1, \nm{x_2} _2,\dots)$.

	Then  $X$  is a \bspc\ and we have inclusions
	$I_n:X_n\to X$,  $I_n(x_n)=(0,\ldots,0,x_n,0\ldots)$ with
	$x_n$ at the $n$th place
	and projections
	$\pi_n:X\to X_n$,  $\pi_n(x)=x_n$,
	$\pi_{m,n}:X\to X$,
	$\pi_{m,n}(x)=(z_i)$, where $z_i=x_i$ if $m\le i\le n$, $z_i=0$
	otherwise; $1\le m\le n\le\infty$, not both $\infty$.
	The  $I_n$  are isometries onto
	their image and $I_n$, $\pi_n$, $\pi_{m,n}$  have operator norm
$1$.

\medskip
\proclaim{\hivI}
{
	For a suitable \seq\ of integers  $n(p)\ge1$,  $p\ge2$, and for
	any $Y$ as above, on the $Y$-sum $X$ of  
	$(\lpee(\CC^{\,n(p)}))_{ p = 2}^\infty$  
	there exists a \delbarclosed\  $f\in C^\infty_{0,1}(B_X(1))$
	for which \dbe is not solvable on any  
	$B_X(a,r)\subset B_X(1)$.
}
\endproclaim
	{\sc Remark.}  For $1\le p<\infty$ regard $\lpee(\CC^{\,n})$
	embedded in $\lpee$ via
	$J_n:\lpee(\CC^{\,n})\to\lpee$, 
	$J_n(x_1,\ldots,x_n)=(x_1,\ldots,x_n,0,\ldots)$,
	put
	$B_p(r)=B_{\lpee}(r)$, $B_{p,n}(r)=B_{\lpee(\CC^{n})}(r)$,
	$\ro_n:\lpee\to\lpee(\CC^{\,n})$, 
	$\ro_n(x)=(x_1,\ldots,x_n,0,\ldots)$,
	and let $f\in C^{\,m}_{0,1}(B_p(1))$ be \delbarclosed\ and of
	finite $C^m(B_p(1))$ norm for some $m\ge1$.
	Then, for some $0<r\le1$, $\delbar u = f$ 
	has a bounded solution $u$ on $B_p(r)$ 
 	if and only if for all $n\ge1$
	there are solutions of $\delbar u_n=J^*_n f$ 
	on $B_{p,n}(r)$
	such that
	$\sup_{B_{p,n}(r)}|u_n|\le M \nm{J^*_n f} _{C^m(B_{p,n}(1))}
	\le M \nm{f} _{C^m(B_p(1))}$, with  $M$
	independent of the dimension $n$.

	This observation is the pillar of the argument below.
	Such a reformulation of the solvability of \dbe was
	already given for Hilbert space by Mazet in [M],
	Appendix~3, Section~1, Remark~2.

\proclaim{\hivII}
{
	With the notations of the remark above,
	the following statement $(E_p)$  is false for any integer $p\ge2$.

	$(E_p)$:  
	There exist 
	a radius  $0<r\le1$,  a constant $0<M<\infty$
	such that for all $n=1,2,\dots$ and for all \delbarclosed\ 
	$f\in\cvne(B_{p,n}(1))$ of finite  $C^{p-1}(B_{p,n}(1))$ norm,
	the equation $\delbar u = f$ has a solution on
	$B_{p,n}(r)$ satisfying 
$$
	\sup_{B_{p,n}(r)} |u| \le M \nm{f} _{C^{p-1}(B_{p,n}(1))} .
 $$
}
\endproclaim
\demo{Proof}	
	Denote by $(E'_p)$ the statement $(E_p)$ with `$f\in\cvne(B_{p,n}(1))$'
	replaced by `$f\in C^{p-1}_{\;0,1}(B_{p,n}(1))$'.
	Fix $p$ and suppose for a contradiction that $(E_p)$ is true.
	Since the $\delbar$ differential operator has constant
	coefficients, approximation using convolution shows that
	$(E'_p)$ is also true.

	We claim that $(E'_p)$  implies the solvability 
	on $B_p(r)$ of \dbe
	with any \delbarclosed\ $f\in C^{p-1}_{\;0,1}(B_p(1))$
	of finite $C^{p-1}(B_p(1))$ norm. Let $u_n$ be a solution of
	$\delbar u_n=f|_{B_{p,n}(r)}$ guaranteed by $(E'_p)$.
	The functions $v_n=\ro^*_n u_n$ on $B_p(r)$ satisfy, with a suitable
	constant $N$, that $|v_n(x)|$, $|(\delbar v_n)(x)\xi|\le N$
	for $x\in B_p(r)$, $\xi\in B_p(1)$.
	
	It follows from the Cauchy--Pompeiu representation formula
	[H, Thm.~1.2.1] applied to 1-dimensional slices that
	$(v_n)_1^\infty$ 
	is a locally equicontinuous family on $B_p(r)$.
	The Arzel\`a--Ascoli theorem gives a subsequence $v_{n'}\to v$
	converging uniformly on compacts in $B_p(r)$.
	As $v$ is continuous and 
	$\delbar v=f$ holds restricted to $B_{p,n}(r)$ for every $n$ in the
        distributional sense, it follows by approximation
	that $\delbar v=f$ holds
	in the distributional sense restricted to any finite dimensional slice
	of $B_p(r)$.
	The ``elliptic regularity of the $\delbar$~operator''
	implies that $v$ is a $C^{p-1}$ solution of $\delbar v =f$
	on $B_p(r)$. See [L2, Props.~2.3, 2.4].
		
	Now, pull back the form  $g$  in Coeur\'e's or Lempert's
	example for  $\lpee$  by  $x\mapsto\epsz x$  with an  $\epsz>0$
	so small that the resulting form  $f$  has finite
	$C^{p-1}(B_{\lpee}(1))$  norm.  Then \dbe is not solvable 
	on any open subset of  $\lpee$. This contradiction proves \hivII.
\enddemo
\demo{Proof of \hivI}  
	We shall use the method of ``Condensation of Singularities.''
	As  $(E_p)$  is false for $p\ge2$ integer, we have
	\seq{}s  $n(p)\ge1$ of integers,  $r_p\to+0$  of radii, 
	$f_p\in\cvne(B_{p,n(p)}(1))$  of \delbarclosed\ forms
	with
	$\nm{f_p} _{C^{p-1}(B_{p,n(p)}(1))}=1$  such that if  $\delbar u=f_p$	
	on  $B_{p,n(p)}(r_p)$  then  $\sup_{B_{p,n(p)}(r_p)}|u|\ge p^{p+1}$.

	Let  $X$  be the  $Y$-sum of  $\lpee(\CC^{\,n(p)})$, $p=2,3,\ldots$.
	Put  $f=\sum_{p=2}^\infty p^{-p} \pi_p^* f_p$.
	One checks that  $f$ is in $\cvne(B_X(1))$ and  is  \delbarclosed.

	We claim that  $\delbar u=f$  cannot be solved on any open subset
	of  $B_X(1)$.
	
	Indeed, suppose for a contradiction that there are a ball  $B_X(a,r)$
	and a function  $u$  with  $\delbar u=f$  on  $B_X(a,r)$.
	Take  $r$  so small that  $\sup_{B_X(a,r)}|u|=N<\infty$.
	This can be done as  $u$  is continuous at  $a$  (even \cinfty).
	Choose $q\ge2$ so large that  $\nm{\pi_{q+1,\infty}(a)} <r/3$.
	Fix  $p>q,N$  so large that  $r_p<r/3$.
	Let  $v(z)=u(\pi_{2,q}(a)+I_p(z))$  for  $z\in B_{p, n(p)}(r_p)$.
	Then  $\delbar v=p^{-p} f_p$  on  $B_{p, n(p)}(r_p)$, so
	$N\ge\sup_{B_{p, n(p)}(r_p)} |v|\ge p^{-p}p^{p+1}=p>N$. 
	This contradiction proves \hivI.	
\enddemo
	Further, we claim that
	$\hbox{dim}_\CC H^{0,1}_{\;\delbar}(B_X(1))=\infty$.
	We group the indices $p$ into pairwise disjoint 
	infinite sets $P_n$, $n\ge1$.
	Then for the $Y$-sum $X^{(n)}$ of 
	$\lpee(\CC^{\,n(p)})$, $p\in P_n$, we have inclusions
	$J_n:X^{(n)}\to X$ and projections $\ro_n:X\to X^{(n)}$ both
	of operator norm 1. Let $g_n\in\cvne(B_{X^{(n)}}(1))$ be 
	a \delbarclosed\ nowhere \delbarexact\ form whose 
	existence is guaranteed by the proof of \hivI.
	Then $f_n=\ro_n^*g_n$ are linearly independent in
	$H^{0,1}_{\;\delbar}(B_X(1))$. Indeed, suppose that 
	$\la_1 f_1+\ldots+\la_n f_n=\delbar u$, $\la_i\in\CC$.
	Then by restricting to $X^{(i)}$ we see that $\la_i g_i$ is
	\delbarexact, hence $\la_i=0$.		

	Should it turn out (as it is yet unknown) that on the
	unit ball $B$ of $l_2$ there are \delbarclosed\ \zoforms\
	of arbitrarily high finite smoothness that are nowhere
	\delbarexact, then the construction in Section~1 with $Y=l_2$
	would
	yield a \delbarclosed\ $f\in\cvne(B)$ which is nowhere
	\delbarexact: a non-solvable $\delbar u=f$ in Hilbert space.


\subhead 2.\ \ Approximation
\endsubhead

	We consider the following kind of approximation  
	in a Banach space $X$.

{\it
	$(A)$:  For any  $0<r<R$,  $\epsz>0$, and $f:B_X(R)\to\CC$
	holomorphic, there exists an entire function  $g:X\to\CC$
	with  $|f-g|<\epsz$  on  $B_X(r)$.
}
\proclaim{\hivIII}
	The statement  $(A)$  holds for the  $l_1$-sum 
	$X$  of any sequence of finite
	dimensional Banach spaces  $(X_n, \nm{\cdot} _n)$.
\endproclaim
	Lempert in [L4] proved  $(A)$  for  $X=l_1$. When this
	manuscript was first written, 
	\hivIII\ was the most general theorem
	proving $(A)$. Later, however, 
	$(A)$ was proved in [L6]
	for any $X$ with a countable unconditional basis, i.e., for
	most classical \bspc{}s.
	It is not
	clear whether all spaces $X$ in \hivIII\ have a countable
	unconditional basis, or even a Schauder basis.

	The proof of \hivIII\ is a modification and extension of Lempert's
	method in [L4].
	Lempert's argument is based on the so-called monomial expansion
	of functions holomorphic on a ball  $\Vert x\Vert<R\le\infty$
	of  $\l_1$  (an analogue of the power series
	expansion on a finite dimensional space), and on the use of a 
	dominating function  $\Delta(q,z)$  defined and continuous on
	$\CC\times B_{l_1}(1)$,  whose role in the estimation of monomial
	series is similar to the role of the geometric series in 
	estimating power series.	
	
	We replace the monomials by so-called multihomogeneous
	functions but use the same dominating function  $\Delta$
	of Lempert.

	{\bf 2.1. Multihomogeneous functions.} 	
	Let  $X$  be the $l_1$-sum of a \seq\ of finite \dim{}al
	\bspc{}s  $(X_n, \nm{\cdot} _n)$.
	For  $\lambda=(\lambda_n)\in\l_\infty$  and  $x\in X$ put
	$\lambda x=(\lambda_1 x_1, \lambda_2 x_2,\dots)\in X$.
	In the rest of this Section  $k$  denotes a multiindex.
	A \mind\ $k=(k_n)$  for us is a 
	sequence of integers  $k_n\ge0$ with  $k_n=0$  for  $n$  large enough.
	The support of  $k$  is the set  $\hbox{\rm supp\ }k=\{n:k_n\not=0\}$.
	We define  $\Vert k\Vert=\sum |k_n|$,  and  $\#k$  as the number of
	elements of the support of  $k$.
	For a \seq\ of \cpx\ numbers $\lambda=(\lambda_n)$  we put  	
	$\lambda^k=\lambda_1^{k_1}\lambda_2^{k_2}\ldots\in\CC$,
	a finite product.
	For a multiindex  $k$,  a holomorphic function
	$\fii:B_X(R)\to\CC$  is called  {\sl $k$-\homog\ }
	if  $\fii(\lambda x)=\lambda^k\fii(x)$  for all  $x\in B_X(R)$,
	$\lambda=(\lambda_n)\in l_\infty$ with $|\lambda_n|=1$.
	
	A \kmh\ function  $\fii$  is a \homog\ polynomial of degree $\nm{k} $
	depending only on those
	finitely many variables  $x_n\in X_n$ for which  $n\in\hbox{\rm supp\ }k$.
	In particular,  $\fii$  extends automatically to an entire function on $X$,
	and  $\fii(\la x)=\la^k\fii(x)$  holds for all $x\in X$ and  $\la\in l_\infty$.
	
	We define the norm  $[\fii]$  of a \kmh\ \fv\  $\fii$  by
	$[\fii]=\sup_{\Vert x\Vert\le 1} |\fii(x)|$.  The set of all \kmh\
	\fvs\  $\fii$  for a fixed  $k$  is a finite dimensional Banach space
	with this norm.

\proclaim{\hivIV}
	For $\fii$ \kmh,
	$|\fii(x)|\le[\fii]|x|^k {\Vert k\Vert}^{\Vert k\Vert}k^{-k}$,
	$x\in X$.
\endproclaim
\demo{Proof} 
	If  $x_i=0$  for some  $i\in\hbox{\rm supp\ }k$, then
	$\fii(x)=0$  as seen from the definition.  So we may suppose that
	$\hbox{\rm supp\ }k=\{1,2,\dots,n\}$  and  
	$x_i\not=0$  for  $1\le i\le n$. Put  
$$
	y=\left({{k_1}\over{\Vert k\Vert}} {{x_1}\over{\nm{x_1} _1}},
	\dots,{{k_n}\over{\Vert k\Vert}} {{x_n}\over{\nm{x_n} _n}},
	0,\dots\right)\in X.
 $$
	Then  $\Vert y\Vert = 1$,  so  
	$\ds [\fii]\ge|\fii(y)|={{k^k}{{\Vert k\Vert}^{-\Vert k\Vert}}}
	{{|x|^{-k}}}|\fii(x)|$,  as claimed.
\enddemo
	{\bf 2.2. The dominating \fv\ of Lempert.}
	This \fv\ is defined by the series
$$	\Delta(q,z)=\sum_k {{{\Vert k\Vert}^{\Vert k\Vert}}\over{k^k}}
	|q|^{\#k}|z^k|
 $$  
	for  $(q,z)\in\CC\times B_{l_1}(1)$.
	See [L2, Section 4].
\proclaim{\hivV}
	(a)  The series for  $\Delta$  converges uniformly on compacts
	in $\CC\times B_{l_1}(1)$.

	(b)  For each  $0<\theta<1$  there is an  $\epsz>0$  such that
	$\Delta$  is bounded on  $B_\CC(\epsz)\times B_{l_1}(\theta)$.
\endproclaim
\demo{Proof}  See [L4, Thm.~2.1].
\enddemo
	We remark that the norm of a monomial  $z^k$  on  $l_1$  is
	$[z^k]={{k^k}{{\Vert k\Vert}^{-\Vert k\Vert}}}$  as a 
	simple calculation
	shows. So  $\Delta(q,z)$  can be written as 
	$\Delta(q,z)=\sum {{|z^k|}{[z^k]}^{-1}}|q|^{\#k}$,
	where we add up normalized monomials with a weight counting the
	number of variables in the monomials.

	{\bf 2.3. Multihomogeneous expansions.}
	Let  $T=(\RR/\ZZ)^\infty=\{t=(t_n):0\le t_n<1\}$  be the 
	infinite \dim{}al torus,
	a compact topological group with the product topology and with
	Haar measure  $dt$  of total mass equal to 1.

	For a \holo\ \fv\  $f:B_X(R)\to\CC$  we define the \mhomog\ expansion of
	$f$  by  $f\sim\sum f_k$, with
	$f_k(x)=\int_{t\in T} f(e^{2\pi i\,t} x)\,e^{-2\pi i\, (k\cdot t)}\,dt$,
	where  $k$  is any \mind,  
	$e^{2\pi i\,t} x=(e^{2\pi i\,t_1} x_1,e^{2\pi i\,t_2} x_2,\dots)$
	and  $(k\cdot t)=\sum k_n t_n$,  a finite sum.
	Then  $f_k$  is defined, holomorphic and \kmh\ on  $B_X(R)$.
	We call  $f_k$  the \kmh\ component of the \fv\ $f$.	

	Let  $S=\{\sigma=(\sigma_n):0\le\sigma_n\to 0\}$ and
	$S_1=\{\sigma\in S:0\le\sigma_n<1\}$, 
	$\sa A=\halmaz{\sa x : x\in A} $ for $A\subset X$, $\sa\in S$
	as in [L4, Section~2].

\proclaim{\hivVI}
	(a)  If  $f:B_X(R)\to\CC$  is a \holo\ function, then
	we have the estimate
	$M(\sigma)=\sup_k [f_k]\sigma^k R^{\Vert k\Vert} <\infty$
	for all  $\sigma\in S_1$.

	(b)  If  $f_k$ is \kmh\ and  $M(\sigma)<\infty$  for all  $\sigma\in S_1$,
	then the series  $g=\sum f_k$  converges uniformly on compact subsets of
	$B_X(R)$,  $g$  is \holo\ on  $B_X(R)$,  and  the \kmh\ component
	$g_k$  of $g$  is equal to  $f_k$.
\endproclaim
\demo{Proof}
	We use the following compactness criterion:	
	A subset  $K\subset X$  is compact if and only if  $K$  is closed,
	bounded,
	and the tail sums  $R_n(x)=\sum_{\nu\ge n} \nm{x_\nu} _\nu\to 0$  
	uniformly on $K$ as $n\to\infty$.

	We outline the proof. If $K$ is compact, then $R_n\to0$
	uniformly on $K$ by Dini's \thm\ on monotone convergence of
	continuous \fvs\ on a compact space to a continuous limit. 
	In the other direction, fix $\epsz>0$. We produce a finite
	covering of $K$ by $\epsz$-balls. Fix $n$ so large that
	$R_n<\epsz/2$ on $K$, and project $K$ onto the space of the
	first $n$ coordinates, this is a bounded set in a finite dimensional
	space, so it has a finite covering by balls $B_X(x_i,\epsz/2)$.
	Now, $B_X(x_i,\epsz)$ cover $K$.

	This criterion implies, in particular, that any compact $K\subset B_X(1)$
	is contained  in
	$\sigma^2 B_X(1)$  for suitable  $\sigma\in S_1$, and all the
	sets $\sa B_X(1)$ for $\sa\in S$ have compact closure.
	The utility of such a criterion was already observed by Ryan [R]
	in a similar context.

	{\it Proof of (a).} The set  $\overline {\sigma B_X(R)}$ being compact,
	$\sup_{\nm{x} <1}|f(\sa Rx)|=M<\infty$. Thus,
	$\fii(x)=f(\sa Rx)$ for $\nm{x} <1$ is bounded by $M$ on $B_X(1)$.
	So is its \kmh\ component $\fii_k(x)=\sa^k R^\nm{k}  f_k(x)$,
	hence $[f_k]\sa^k R^{\nm k } \le M$, or $M(\sa)\le M<\infty$.

	{\it Proof of (b).}  \WLOG\ we may suppose that the given
	compact is 
	$\sigma L$,  where  $L\subset B_X(r)$  is compact,
	$\sigma\in S_1$  and  $r<R$.
	Then putting  $x=\sigma y$  for  $|y|<r$,  $y\in L$,  we have that
$$	|f_k(x)|\le[f_k]{\Vert k\Vert}^{\Vert k\Vert} k^{-k}
	|x|^k=[f_k]{\Vert k\Vert}^{\Vert k\Vert} k^{-k} \sigma^k |y|^k=
	[f_k]\sigma^k R^{\Vert k\Vert} \cdot {\Vert k\Vert}^{\Vert k\Vert} k^{-k}
	|y/R|^k.
 $$
	
	Summing on  $k$,  we get
	$\sum |f_k(x)|\le M(\sigma)\Delta(1, z)\le M<\infty$
	where  $z=|y/R|$  ranges in a compact subset of  $B_{l_1}(1)$,  and
	the series for $\Delta$  converges uniformly  by \hivV(a).

	This concludes the proof of \hivVI. 
\enddemo
\proclaim{\hivVII}
	Let $f_k$ be \kmh. 
	If for each \mind\  $k$  and  for all  $\sigma\in S$ (!)
	we have  $\sup_k [f_k]\sigma^k R^{\Vert k\Vert}<\infty$,
	then  $\sum f_k$  is an entire function on  $X$.
\endproclaim
\demo{Proof}  
	If  $M(\sa)<\infty$  for all $\sa\in S$, then
	$M(\la\sa)<\infty$  for all $0<\la<\infty$, $\sa\in S_1$,
	which has the same effect as changing  $R$  to  $\la R$.
	Hence the \mhomog\ series converges on the whole of $X$.
\enddemo
	We quote two lemmas from [L4].
\proclaim{\hivVIII}
	If the numbers  $0\le c_k<\infty$  are such that  
	$\sup_k c_k\sigma^k<\infty$  for all  $\sigma\in S_1$,  then for any
	$Q\ge 1$ and $\sigma\in S_1$  the estimate 
	$\sup_k c_k\sigma^kQ^{\#k}<\infty$  holds.
\endproclaim
\demo{Proof}  See [L4, Prop.~4.2].
\enddemo
\proclaim{\hivIX}
	Let  $0<\theta<1$  and  $\KK$  a set of multiindices  $k$.
	Then if  $0<c_k<\infty$, $k\in\KK$, satisfy  
	$\inf_{k\in\KK} c_k \theta^{\Vert k\Vert}>0$  and
	$\sup_{k\in\KK}	c_k \sigma^k<\infty$  for all  $\sigma\in S_1$,
	then  $\sup_{k\in\KK} c_k \sigma^k<\infty$  for all  $\sigma\in S$, too.
\endproclaim
\demo{Proof}  See [L4, Prop.~4.3].
\enddemo
\demo{Proof of \hivIII} 
	Let us expand  $f$  in  a  \mhomog\ series  $\sum f_k$.
	Fix  a number $0<\theta<1$  with  $r<\theta^2 R$.
	For any  $\delta>0$,  $Q>1$  (to be suitably chosen below)  put
	$c_k=[f_k]R^{\Vert k\Vert}$,  $c'_k=c_k Q^{\#k}$,
	$\KK=\left\{k:c'_k\theta^{\Vert k\Vert}\equiv
	[f_k](\theta R)^{\Vert k\Vert} Q^{\#k}\ge\delta\right\}$, and
	$g(x)=\sum_{k\in\KK} f_k(x)$.

	We claim that this  $g$  is an entire function on  $X$.

	Indeed, by \hivVII\ it is enough to show for all  $\sigma\in S$
	that
$$	\sup_{k\in\KK} [f_k]\sigma^k R^{\Vert k\Vert}\equiv
	\sup_{k\in\KK} c_k \sigma^k<\infty.
 $$
	As $\inf_{k\in\KK} c'_k \theta^{\Vert k\Vert}\ge\delta>0$, and for
	$\sigma\in S_1$  \hivVI(a) implies that
$$	\sup_{k\in\KK} [f_k]\sigma^k R^{\Vert k\Vert}\equiv
	\sup_{k\in\KK} c_k\sigma^k<\infty,
 $$  
	so by \hivVIII, 
	$\sup_{k\in\KK} c'_k\sigma^k<\infty$  holds for all  $\sigma\in S_1$.
	Now both conditions of \hivIX\ are verified, hence
	$\sup_{k\in\KK} c_k\sigma^k \le \sup_{k\in\KK} c'_k \sigma^k<\infty$  
	for all  $\sigma\in S$. Therefore, by \hivVII,
	 $g$  is an entire function.

	For  $k\not\in\KK$ we have  $[f_k](\theta R)^{\Vert k\Vert} Q^{\#k}\le\delta$.
	We estimate  $|f(x)-g(x)|$.  For  $\Vert x\Vert<r$ by \hivIV, we have
$$
\eqalign{
	|f(x)-g(x)|&\le\sum_{k\not\in\KK} |f_k(x)| \le \sum_{k\not\in\KK} [f_k]
	{{{\Vert k\Vert}^{\Vert k\Vert}}\over{k^k}}|x|^k \cr
        &\le\sum_{k\not\in\KK} \delta Q^{-\#k} (\theta R)^{-\Vert k\Vert}
	{{{\Vert k\Vert}^{\Vert k\Vert}}\over{k^k}}|x|^k \cr &\le
	\delta\sum_{k\not\in\KK} Q^{-\#k} {{{\Vert k\Vert}^{\Vert k\Vert}}\over{k^k}}
	\left|{{\theta x}\over r}\right|^k  
	\le\delta\sup_{\Vert w\Vert\le\theta} \Delta(Q^{-1},w)
}
 $$
	as $\theta R>r/\theta$ and  $|w|=|\theta x/r|\le\theta$.  
	But the last expression can be made $<\epsz$
	by 
	choosing first $Q$  large enough to make the  $\sup$  finite by \hivV(b),
	and then by choosing  $\delta$  small enough.
	
	Thus, the proof of the approximation \hivIII\ is concluded. 
\enddemo
	Let $Y=l_q$, $1\le q<\infty$, or $Y=c_0$. 	
	Let $e_{pi}^n$, $1\le i\le n$ be the standard basis of
	$\lpee(\CC^{\,n})$.
	Then the $Y$-sum $X$ of any sequence $\l_{p_k}(\CC^{\,n_k})$
	spaces, $k\ge1$,
	has a countable unconditional basis:
	$e_{p_1 1}^{n_1}$, $e_{p_1 2}^{n_1}$, $\ldots$, $e_{p_1 n_1}^{n_1}$;
	$e_{p_2 1}^{n_2}$, $\ldots$, $e_{p_2 n_2}^{n_2}$; $\ldots$.
	Now, the approximation theorem of Lempert [L6, Thm.~0.1], or 
	in the case $Y=l_1$, \hivIII\ above, implies by the
	vanishing theorem [L3, Thm.~0.3]
	that the sheaf cohomology groups 
	$H^q(\Omega, \OO)=0$, $q\ge1$, on any pseudoconvex open set
	$\Omega\subset X$ for the sheaf $\OO$ of germs of \holo\
	\fv{}s on $X$.  So for any $Y$, the space $X$ of \hivI\ has the
	property that $H^{0,1}_{\;\delbar}(\Omega)\not=0$ (in fact, 
	infinite dimensional) and $H^1(\Omega, \OO)=0$ for any
	bounded pseudoconvex open set $\Omega\subset X$: the Dolbeault
	\isom\ \thm\ does not
        generalize to arbitrary \bspc{}s.

	We remark that if the form $f$ is real-analytic and $\Omega$
	pseudoconvex, then by [L1, Prop.~3.2] the equation \dbe has 
	real-analytic local 
	solutions; since $H^1(\Omega, \OO)=0$, we get a 
	global real-analytic solution, too.


\subhead 3.\ \ \Acc\ \mfold{}s
\endsubhead

	\hivI\ verifies 
	the hypothesis of \hivX\ below in a case.

\proclaim{\hivX}
	Let  $X$  be a \bspc\ and suppose that on  $B=B_X(1)$
	there exists a \delbarclosed\   $f\in\cvne(B)$  that is
	not \delbarexact\ on any open subset.  Then on  $M=B\times\CC$
	a \cinfty-smooth integrable \acc\ \struc\ 
	$M_f$
	can be constructed 
	in such a way that no open subset of  $M_f$  is bi\holo\
	to an open subset of a \bspc.
\endproclaim

	As the referee  	
	has kindly pointed it out, 
	the method of this section is
	analogous to one used earlier to
	construct nonrealizable CR hypersurfaces by Jacobowitz in~[J].
	
	We recall the definition of \acc\ \struc.
	An \acc\ \struc\ on a  $C^m$-smooth manifold  $M$  is a
	splitting of the complexified tangent bundle 
	$\ctm=T^\en\oplus T^\ne$  into the direct sum of
	two complex vector bundles of class  $C^{m-1}$
	with  $T^\ne=\overline{T^\en}$, $m=1,\ldots$,~$\infty$,~$\omega$ and
	$m-1=m$ for $m=\infty$,~$\omega$.
	An \acc\ \struc\  is called formally integrable (or just integrable)
	if  $m\ge2$ and the Lie bracket of two $(1,0)$ vector fields of 
	class $C^1$ is also a $(1,0)$ vector field; here $(1,0)$ can
	be changed to $(0,1)$.

	The proof of \hivX\ requires a few steps.

	{\bf 3.1. Construction of the \acc\ \struc\ on  $M$.}
	The construction will be described in a setting more general
	than that of \hivX, namely, in the context of principal
	bundles.

	Denote by $\za^\en$, $\za^\ne$ the $(1,0)$-part, $(0,1)$-part
	of a complex tangent vector $\za$ to an \acc\ \mfold.
	Let $B$ be a complex Banach manifold, $G$ a 
	finite dimensional \cpx\ Lie \gp\
	with Lie algebra $\gg=T_eG$, $f\in\cvne(B,\gg^\en)$ a
	\zoform\ with values in $\gg^\en$, and $L_z:G\to G$ the 
	left translation $L_z(s)=zs$, $z,s\in G$.
	Define the \holo\ Maurer--Cartan form 
	$\mu\in{C^{\,\infty}_{1,0}}(G,\gg^\en)$ 
	by 
$$	
	\mu(\nu)=(dL_z)^{-1}\nu^\en=((dL_z)^{-1}\nu)^\en
 $$
	for $\nu\in\CC\otimes T_z G$.
	Recall the \holo\ Maurer--Cartan formula
	$d\mu+\ek[\mu,\mu]=0$, which can be proved similarly to
	or deduced from the usual Maurer--Cartan formula.
	Define on any complex Banach manifold $N$ the Lie bracket
	$[\fii,\psi]\in C^\infty_{0,2}(N,\gg^\en)$
	of forms  
	$\fii,\psi\in\cvne(N,\gg^\en)$ 
	by the usual formula 
$$
	[\fii,\psi](\za,\za')=[\fii(\za),\psi(\za')]-[\fii(\za'),\psi(\za)]
 $$
	for $\za,\za'\in \CC\otimes T_x N$,
	where the brackets on the right hand side are taken in the
	Lie algebra $\gg^\en$.
	In particular, $[f,f](\za,\za')=2\,[f(\za),f(\za')]$.

	We define an \acc\ \struc\ $M_f$ on $M=B\times G$ by putting 
	$(\za,\nu)\in \CC\otimes T_{(x,z)}M=
	\CC\otimes T_x B\oplus\CC\otimes T_z G$ 
	in $T^{\;\ne}_{(x,z)}M$
	if and only if 
$$
	\za=\za^\ne\qquad \hbox{and} \qquad\mu(\nu)=f(\za).\tag\acm
 $$

	In the setting of \hivX\ we identify $G=\CC$ and $\gg^\en=\CC$ via
	the correspondence $G\ni s \sim 
	s\, \partial/(\partial z)|_{z=0}\in\gg^\en$,
	where $z$ is
	the usual coordinate on $\CC$.

	We verify below the following:
	Definition (\acm)  gives an \acc\ \struc\ $M_f$
	on $M$ and makes it into an \acc\ principal $G$ bundle;
	$M$ is formally integrable if and only if $\delbar f+\ek\,[f,f]=0$
	holds; if $M_f$ is locally bi\holo\ to a \bspc\ then 
	$\Dbar u=f$ where $u:B\to G$ is defined locally
	and $\Dbar$ is defined by 
$$
	\Dbar u\;(\za)=\mu\left(du\;(\za^\ne)\right)
 $$
	for $\za\in T_x B$.
	In the setting of \hivX\ this $\Dbar u$ reduces to the usual $\delbar u$.

	{\bf 3.2. Verification.}  To verify that (\acm) defines an
	almost complex structure on  $M$, we need to check
	conditions $1^\circ$--$2^\circ$.

	$1^\circ$ {\sl If $V=(\za,\nu)$, 
	$\Vbar=(\bbar{\za},\bbar{\nu})$ are in $T^{\;\ne}_{(x,z)}M$
	then $V=0$.
	}
	We have $0=\za^\ne=\bbar{\za}^\ne\equiv\overline{\za^\en}$, or $\za=0$.
	Similarly $\mu(\nu)=\mu(\bbar{\nu})=0$ implies
	$\nu^\en=\bbar{\nu}^\en=0$,
	or $\nu=0$.

	$2^\circ$ 
{\sl
	Given $V=(\za,\nu)$, decompose it as $V=V_1+V_2$ with 
	$V_1, \Vbar_{\!2}\in T^\ne M$.
}
	One checks that 
$$
\eqalign{
	V_1=\left(\za^\ne,\; dL_z\;f(\za) - dL_z\;\overline{f(\bbar{\za})}
				 + \nu^\ne\right) \cr
	V_2=\left(\za^\en,\; dL_z\;\overline{f(\bbar{\za})} - dL_z\;f(\za)
                                 + \nu^\en\right) 
}
 $$	
	is the unique way of decomposition.

	$3^\circ$ 
	Condition of formal integrability: 
{\sl 
	If $V=(\za,\nu)$, $V'=(\za',\nu')$ are \cinfty\ sections 
	of $T^\ne M$ over an open subset of
	$M$, then their Lie bracket $[V, V']$ is also a section of $T^\ne M$.
}

	Denote by $\LL_Z$ the Lie derivative along a complex
	vector field $Z$. 
	We can write $[V,V']$ as
$$
 	[V,V']=(\za^*,\nu^*)=
	\left([\za,\za']+\LL_\nu\za'-\LL_{\nu'}\za,\;
	      [\nu,\nu']+\LL_\za\nu'-\LL_{\za'}\nu\right).
 $$

	We work out below the condition of formal integrability for $M_f$
	in terms of $f$.

	(a) The first component $\za^*$ is $(\ne)$  because so are
	$[\za,\za']$, $\LL_\nu\za'$, $\LL_{\nu'}\za$ since $B$ is 
	a \cpx\ manifold.
	
 	(b) Taking the $\LL_\za$, $\LL_\nu$ Lie derivatives of the identity
	$\mu(\nu')-f(\za')=0$ and reversing the roles of $V,V'$ we find 
	the equations
$$
\alignat{3}
	& \mu(\LL_\za\nu') &-& \LL_\za(f(\za')) &=& 0 \\
	& \mu(\LL_{\za'}\nu) &-& \LL_{\za'}(f(\za)) &=& 0 \\
	& \LL_\nu(\mu\nu') &-& f(\LL_\nu\za')   &=& 0 \\
	& \LL_{\nu'}(\mu\nu) &-& f(\LL_{\nu'}\za)   &=& 0 \\
\endalignat
 $$
	whose alternating sum is
$$
\eqalign{
	\{\mu\nu^*-f(\za^*)\} & + 
	\{\LL_\nu(\mu\nu')-\LL_{\nu'}(\mu\nu)-\mu([\nu,\nu'])\} - {} \cr
	& \{\LL_\za(f(\za'))-\LL_{\za'}(f(\za))-f([\za,\za'])\} = 0.
}
 $$
	Hence, by Cartan's formula 
	for exterior derivatives,
	the condition of formal integrability is that 
$$
	(d\mu)(\nu,\nu')-(df)(\za,\za')=0.
 $$
	Since
	$\mu(\nu)=f(\za)$, $\mu(\nu')=f(\za')$ we get 
	by the  \holo\ Maurer--Cartan formula
	that
$$
	-\ek\,[f,f](\za,\za')-(df)(\za,\za')=0
 $$
	for all $(0,1)$ vector fields $\za$, $\za'$ on $B$.
	Hence the \acc\ manifold $M_f$ is formally integrable if and only if
$$
	\delbar f + \ek\,[f,f] = 0,
 $$
	which condition reduces to $\delbar f=0$ when $G$ is commutative
	as in \hivX\.


	{\bf 3.3. Geometric properties of $M$.} 
	To check that $M$ is
	a principal $G$ bundle we need to
	verify that $\pi:M=B\times G\to B$, $\pi(x,z)=x$ is \holo\
	and that $G$ has a simply transitive action on the fibers of $M$.
	Indeed, $\pi$ is \holo\ as $d\pi(\za,\nu)=\za$
	takes $(0,1)$-vectors to $(0,1)$-vectors.
	The action of $w\in G$ on $M$ is given by
	the left translation $l_w(x,z)=(x,wz)$ in the fiber direction. This is
	\holo\ because 
	$(dl_w)(\za,\nu)=(\za,dL_w\;\nu)$ and
	$\mu(dL_w\;\nu)=\mu(\nu)$.

	In the setting of \hivX\ a direct verification shows that
	$\Phi:M_f\to M_g$, $\Phi(x,z)=(x,z+u(x))$ is a bundle biholomorphism,
	where $g=f+\delbar u$ and $u\in C^\infty(B)$ is any function.
	Hence the bundle biholomorphism type of $M_f$ depends only on the
	Dolbeault cohomology class of $f$.

	We return now to the general setting.

\proclaim{\hivXI}
	If $p_0=(x_0,z_0)\in M_f$ has a \nbd\ that is $C^m$-bi\holo,
	$m=1,2,\ldots,\infty$, to an open set in a \bspc,
	then there are a \nbd\ $U_0\subset B$ of $x_0$ and 
	$u\in C^m(U_0,G)$ such that $\Dbar u=f$ on $U_0$.
\endproclaim

\demo{Proof}
	The \bspc\ $\Tne=T_{\;p_0}^\ne M$ has a natural splitting as
	a direct sum $\Tne=\Vne\oplus \Hne$ of 
	vertical and horizontal closed subspaces
$$
\eqalign{
	\Vne = \halmaz{ (0,\nu)\in\Tne : \nu\in T^\ne_{\;z_0}G } \cr
	\Hne = \halmaz{ (\za,\nu)\in\Tne : \nu\in T^\en_{\;z_0}G }
}.
 $$

	Suppose now that $\Phi:U\to V$ is bi\-holo\-morphism of a 
	\nbd\ $U$ of $p_0$ in $M_f$ onto a \nbd\ $V$ of $0$ in a \bspc\ $W$.
	Then the splitting $\Tne=\Vne\oplus \Hne$ induces via 
	$(d\Phi)(p_0)$ a splitting $T^\ne_{0}W\equiv W = \VWne\oplus \HWne$.
	Since $N=\Phi^{-1}(V\cap\HWne)$ is an \acc\ $C^m$-submanifold
	of $M$ passing through $p_0$ transversely to $\Vne$, 
	hence to $\halmaz{x_0} \times G$,
	$N$ is the image near $p_0$ of a \holo\ section $s:U_0\to G$
	on a \nbd\ $U_0$ of $x_0$. Then writing $s(x)=(x,u(x))$ and applying
	(\acm) we obtain that $\mu(du\,(\za))=f(\za)$ for all
	$\za\in T^\ne_{\;x}U_0$, but this is the same as saying
	$\Dbar u=f$ on $U_0$; thus concluding the proof of \hivXI\ and hence 
	that of \hivX.
\enddemo

	We have seen that the Newlander--Nirenberg \thm\ does not generalize to
	arbitrary integrable \acc\ Banach \mfold{}s. It is unknown if it
	generalizes to Hilbert manifolds.

	{\sc Acknowledgement.} The author is grateful to his advisor,
	L.~Lempert for his generous help and kind encouragement.




\enddocument

\bye